\documentclass{tac}
\usepackage{amsmath,latexsym,amssymb,mathrsfs}
\usepackage{textcomp}
\usepackage{xypic}
\usepackage[matrix]{xy}
\input diagxy

\makeatletter
\renewcommand*\env@matrix[1][*\c@MaxMatrixCols c]{%
  \hskip -\arraycolsep
  \let\@ifnextchar\new@ifnextchar
  \array{#1}}
\makeatother

\author {Thomas Weighill}
\address{Centre for Mathematical Structures \\ Department of Mathematical Sciences \\ Stellenbosch University \\ Private Bag X1 \\ 7602 Matieland \\ South Africa \\ $\;$  \\ $\;$ \\  Department of Mathematics \\ The University of Tennessee \\ 227 Ayres Hall \\ 1403 Circle Drive \\ Knoxville, TN 37996-1320 \\ United States of America}
\title{Mal'tsev objects, $R_1$-spaces and ultrametric spaces}
\copyrightyear{2015}
\eaddress{tweighil@vols.utk.edu}

\newtheorem{lemma}{Lemma}

\newtheorem{question}{Question}
\mathrmdef[Hom]{hom}
\mathbfdef{Ord}
\mathbfdef{Top}

\newcommand{\columnII}[2]{\fontsize{9}{10}\selectfont \begin{pmatrix} #1 \\ #2 \end{pmatrix}}

\newcommand{\rowII}[2]{\fontsize{9}{10}\selectfont \begin{pmatrix} #1,  #2 \end{pmatrix}}
\newcommand{\rowIII}[3]{\fontsize{9}{10}\selectfont \begin{pmatrix} #1,  #2,  #3\end{pmatrix}}
\newcommand{\matrIIxIII}[6]{\fontsize{9}{10}\selectfont \begin{pmatrix} #1 & #4 \\ #2 & #5 \\ #3 & #6\end{pmatrix}}
\newcommand{\matrIIxII}[4]{\fontsize{9}{10}\selectfont \begin{pmatrix} #1 & #3 \\ #2 & #4 \end{pmatrix}}
\newcommand{\matrIIIxII}[6]{\fontsize{9}{10}\selectfont \begin{pmatrix} #1 & #2 & #3 \\ #4 & #5 & #6\end{pmatrix}}

\begin{document}
\maketitle
\begin{abstract}
In this paper we introduce a notion of Mal'tsev object, and the dual notion of co-Mal'tsev object, in a general category. In particular, a category $\mathbb{C}$ is a Mal'tsev category if and only if every object in $\mathbb{C}$ is a Mal'tsev object. We show that for a well-powered regular category $\mathbb{C}$ which admits coproducts, the full subcategory of Mal'tsev objects is coreflective in $\mathbb{C}$. We show that the co-Mal'tsev objects in the category of topological spaces and continuous maps are precisely the $R_1$-spaces, and that the co-Mal'tsev objects in the category of metric spaces and short maps are precisely the ultrametric spaces. 
\end{abstract}

\section{Introduction}
A variety $\mathbb{X}$ of universal algebras is called a \emph{Mal'tsev variety}~\cite{Smith} if it satisfies the following condition:

\begin{itemize}
\item[$(\mathsf{M}_1)$] the algebraic theory of $\mathbb{X}$ contains a ternary term $\mu$ satisfying the term equations
$$
\mu(x,y,y) = x = \mu(y,y,x).
$$
\end{itemize} 
A famous theorem of Mal'tsev states that these varieties are precisely those in which the composition of congruences on any object is commutative~\cite{Mal54}. The notion of \emph{Mal'tsev category} is a generalisation of the notion of Mal'tsev variety. Recall that a Mal'tsev category was originally defined in~\cite{Carboni1990} to be a category $\mathbb{C}$ which is exact in the sense of Barr~\cite{Barr1971} and which satisfies the following condition:

\begin{itemize}
\item[$(\mathsf{M}_2)$] every reflexive internal relation in $\mathbb{C}$ is an equivalence relation.
\end{itemize}
In the present paper, by a Mal'tsev category, we mean (as in~\cite{BoJan08}) a category $\mathbb{C}$  which satisfies the following relational reformulation of $(\mathsf{M}_1)$ due to Lambek~\cite{Lambek}:

\begin{itemize}
\item[$(\mathsf{M}_3)$] every internal relation $R$ in $\mathbb{C}$ is \emph{difunctional}, i.e. it satisfies 
$$(x_1 R y_2 \wedge x_2 R y_2 \wedge x_2 R y_1) \Rightarrow x_1 R y_1.$$
\end{itemize}
Note that conditions $(\mathsf{M}_2)$ and $(\mathsf{M}_3)$ can both be formulated in a general category. Recall that an internal relation from an object $X$ to an object $Y$ in a category $\mathbb{C}$ is a triple $(R, r_1, r_2)$ with $R$ an object of $\mathbb{C}$ and $r_1: R \rightarrow X$ and $r_2: R \rightarrow Y$ morphisms of $\mathbb{C}$ such that $r_1$ and $r_2$ are jointly monomorphic. Note that if $\mathbb{C}$ admits binary products, then an internal relation from $X$ to $Y$ can also be viewed as a monomorphism $(r_1, r_2): R \rightarrow X \times Y$. We say that a relation $(R, r_1, r_2)$ from $X$ to $Y$ is reflexive, symmetric, transitive or difunctional when for every object $S$ of $\mathbb{C}$, the relation 
$$ \bfig
\morphism(0,0)<-1000, 0>[\Hom(S, R)`\Hom(S, X) ;\Hom(S, r_1)] 
\morphism(0,0)<1000, 0>[\Hom(S, R)`\Hom(S, Y);\Hom(S, r_2)]
\efig $$
between sets is reflexive, symmetric, transitive or difunctional in the usual sense. For a category $\mathbb{C}$ with finite limits, $(\mathsf{M}_2)$ is equivalent to $(\mathsf{M}_3)$ (see~\cite{Carboni92}).

A category $\mathbb{C}$ is thus a Mal'tsev category if and only if every object $S$ in $\mathbb{C}$ satisfies the following condition:
\begin{itemize}
\item[$(\mathsf{D})$] for any internal relation $(R, r_1, r_2)$ from an object $X$ to an object $Y$, the following relation is difunctional: 
$$ \bfig
\morphism(0,0)<-1000, 0>[\Hom(S, R)`\Hom(S, X) ;\Hom(S, r_1)] 
\morphism(0,0)<1000, 0>[\Hom(S, R)`\Hom(S, Y);\Hom(S, r_2)]
\efig $$
\end{itemize}
For a general category $\mathbb{C}$, we will call an object $S$ satisfying $(\mathsf{D})$ above a \emph{Mal'tsev object}. Note that the Mal'tsev objects in $\mathbb{C}$ are precisely those objects $S$ for which the functor $\Hom(S, -)$ is $M$-closed in the sense of~\cite{JanClosed}, where $M$ is the matrix
$$
M =  \begin{pmatrix}[ccc|c]
  y & x  & x & y \\
  u & u & v & v \\
  \end{pmatrix}.
$$
We call an object $S$ in $\mathbb{C}$ a \emph{co-Mal'tsev object} if it is a Mal'tsev object as an object of the dual category $\mathbb{C}^\mathsf{op}$. We denote the full subcategory of Mal'tsev objects in $\mathbb{C}$ by $\mathsf{Mal}(\mathbb{C})$. 

In this paper, we first give a characterisation of Mal'tsev objects in the case when $\mathbb{C}$ is a category satisfying certain conditions. This characterisation is based on recent work by Bourn and Z. Janelidze~\cite{BoJan08}. We then show that, for a regular category $\mathbb{C}$ with binary coproducts, $\mathsf{Mal}(\mathbb{C})$ contains every full subcategory of $\mathbb{C}$ which is a Mal'tsev category and which is closed under binary coproducts and regular quotients in $\mathbb{C}$. In Section~\ref{r1}, we show that the co-Mal'tsev objects in $\mathbf{Top}$ (the category of topological spaces and continuous maps) are precisely the $R_1$ spaces~\cite{Davis1961}, i.e.~topological spaces satisfying the following ``separation axiom'':
\begin{itemize}
\item[($\mathsf{R}_1$)] for all $x,y \in X$, if there exists an open set $A$ such that $x \in A$ and $y \notin A$, then there exist disjoint open sets $B$ and $C$ such that $x \in B$ and $y \in C$.
\end{itemize}
In Section~\ref{ultra}, we consider the category $\mathbf{Met}$ of metric spaces and short maps, and show that the co-Mal'tsev objects in this category are precisely the ultrametric spaces, i.e. metric spaces $X$ satisfying
$$
d(x,z) \leqslant \mathsf{max}\{(d(x,y), d(y,z)\}
$$
for any $x,y,z$ in $X$. A classical example of an ultrametric space is the set of rationals $\mathbb{Q}$ equipped with the metric arising from the $p$-adic norm for some prime $p$.

\section{General properties of Mal'tsev objects}\label{general}
By a \emph{(regular) quotient/subobject} of an object $X$ in a category $\mathbb{C}$ we mean a (regular) epi/mono with domain/codomain $X$. We say that a subcategory $\mathbb{D}$ in $\mathbb{C}$ is \emph{closed under regular quotients/subobjects} in $\mathbb{C}$ if the codomain/domain of every regular quotient/subobject of an object in $\mathbb{D}$ is in $\mathbb{D}$. 

\begin{proposition} \label{colimits}
For any category $\mathbb{C}$, $\mathsf{Mal}(\mathbb{C})$ is closed under colimits and regular quotients in $\mathbb{C}$. 
\end{proposition}
\begin{proof}
This can be deduced from general considerations via the Yoneda embedding, but it is also easy to prove directly, as we now show. Let $D: \mathbb{D} \rightarrow \mathbb{C}$ be a diagram whose image is contained in $\mathsf{Mal}(\mathbb{C})$ and which has a colimit $C$ in $\mathbb{C}$, and let $(R, r_1, r_2)$ be an internal relation from $X$ to $Y$. Let $x_1, x_2: C \rightarrow X$ and $y_1, y_2: C \rightarrow Y$ be morphisms such that $x_1 R y_2$, $x_2 R y_2$ and $x_2 R y_1$. Then for each object $A$ in $D$, $x_1 \iota_A R y_2 \iota_A$, $x_2 \iota_A R y_2 \iota_A$ and $x_2 \iota_A R y_1 \iota_A$, where $\iota_A$ is the colimit injection. Since each $D(A)$ is a Mal'tsev object, we have $x_1 \iota_A R y_1 \iota_A$ for each object $A$ in $\mathbb{D}$. Thus there is a family of morphisms $h_A: D(A) \rightarrow R$ such that $r_1 h_A = x_1 \iota_A$ and $r_2 h_A = y_1 \iota_A$ for every object $A$ in $\mathbb{D}$. Using the fact that $r_1$ and $r_2$ are jointly monic, it follows that the morphisms $h_A$ induce a morphism $h: C \rightarrow R$, and it is easy to check that $r_1h = x_1$ and $r_2 h = y_1$, so $x_1R y_1$ as required. Suppose now that $S$ is a Mal'tsev object and that $q: S \rightarrow T$ is a regular epimorphism, which is the coequalizer of $a,b: Q \rightarrow S$. Let $u_1, u_2: T \rightarrow X$ and $v_1, v_2: T \rightarrow Y$ be morphisms such that $u_1 R v_2$, $u_2 R v_2$ and $u_2 R v_1$. Then $u_1q R v_2q$, $u_2q R v_2q$ and $u_2q R v_1q$. Since $S$ is a Mal'tsev object, we have that $u_1q R v_1q$. In other words, there is a map $f: S \rightarrow R$ such that $r_1 f = u_1q$ and $r_2 f = v_1q$. But then $r_1fa = u_1qa = u_1qb = r_1fb$ and $r_2fa = v_1qa = v_1qb = r_2fb$, so since $r_1$ and $r_2$ are jointly monic, we have $fb = fa$. Thus there is a morphism $g: T \rightarrow R$ such that $gq = f$, and one checks that $r_1g = u_1$ and $r_2 g = v_1$, which gives $u_1 R v_1$ as required. 
\end{proof}

Recall that a category $\mathbb{C}$ is \emph{well-powered} if for every object $X$ of $\mathbb{C}$, the collection of all isomorphism classes of subobjects of $X$ may be labelled by a set.

\begin{corollary}\label{reflect}
Let $\mathbb{C}$ be a well-powered regular category which admits coproducts. Then $\mathsf{Mal}(\mathbb{C})$ is a coreflective subcategory of $\mathbb{C}$. In particular, $\mathsf{Mal}(\mathbb{C})$ will be (finitely) complete if $\mathbb{C}$ is (finitely) complete. 
\end{corollary}
\begin{proof}
This follows from the following general fact: if $\mathbb{B}$ is a full subcategory of a well-powered regular category $\mathbb{X}$ which admits coproducts and $\mathbb{B}$ is closed under coproducts and regular quotients, then $\mathbb{B}$ is a coreflective subcategory. Indeed, if $X$ is any object in $\mathbb{X}$, take a set of representatives $M$ of all subobjects of $X$ which lie in $\mathbb{B}$, and let $\coprod M$ be the coproduct of their domains. The coreflection of $X$ into the subcategory $\mathbb{B}$ is then given by the domain of the mono part of the factorisation of the canonical morphism from $\coprod M$ to $X$.
\end{proof}

Proposition~\ref{bojan1} below follows straightforwardly from the proofs of Proposition 4.1 and Theorem 4.2 in~\cite{BoJan08}, but we present a direct proof here for the sake of completeness. For convenience, given an internal relation $(R, r_1, r_2)$ from $X$ to $Y$ and two morphisms $f: S \rightarrow X$ and $g: S \rightarrow Y$, we write $fRg$ to mean that $f$ and $g$ are related by the image of $(R, r_1, r_2)$ under $\hom(S, -)$. Dually, given an internal co-relation $(R, r_1, r_2)$ from $X$ to $Y$ (i.e~a pair of jointly epimorphic morphisms $r_1: X \rightarrow R$ and $r_2: Y \rightarrow R$) and two morphisms $f: X \rightarrow S$ and $g: Y \rightarrow S$, we write $fRg$ to mean that $f$ and $g$ are related by the relation
$$ \bfig
\morphism(0,0)<-1000, 0>[\Hom(R, S)`\Hom(X, S) ;\Hom(r_1, S)] 
\morphism(0,0)<1000, 0>[\Hom(R, S)`\Hom(Y, S);\Hom(r_2, S)]
\efig $$

\begin{proposition}\label{bojan1}
Let $\mathbb{C}$ be a regular category which admits binary coproducts. Then for any object $S$ in $\mathbb{C}$, the following are equivalent:
\begin{itemize}
\item[(a)] $S$ is a Mal'tsev object;
\item[(b)] $\iota_1 R' \iota_1$, where $\iota_1: S \rightarrow 2S$ is the first coproduct injection and $(R', r'_1, r'_2)$ is the internal relation from $2S$ to $2S$ appearing in the (regular epi, mono)-factorisation of the vertical morphism in the following diagram:
\begin{equation}\label{eq:fact1}
\bfig
\Dtriangle(0,0)|aab|<400,400>[3S`R'`2S\times 2S;%
\matrIIxIII{\iota_1}{\iota_2}{\iota_2}{\iota_2}{\iota_2}{\iota_1}`e`r'=(r'_1,r'_2)]
\efig 
\end{equation}
\end{itemize}
\end{proposition}
\begin{proof}
(a)$\Rightarrow$(b): For the internal relation $R'$ in diagram (\ref{eq:fact1}) we have $\iota_1 R' \iota_2$, $\iota_2 R' \iota_2$ and $\iota_2 R' \iota_1$, so $\iota_1 R' \iota_1$ by difunctionality of $\Hom(S,R')$.

(b)$\Rightarrow$(a): Suppose $(R, r_1, r_2)$ is an internal relation from $X$ to $Y$ and $x_1, x_2: S \rightarrow A$ and $y_1, y_2: S \rightarrow B$ are morphisms such that $x_1 R y_2$, $x_2 R y_2$ and $x_2 R y_1$. Consider the diagram of solid arrows:
$$\bfig
\square(0,0)|alrb|/{-->}`{->}`{->}`{->}/<1000,500>[3S`R`2S \times 2S`X\times Y;%
p`\matrIIxIII{\iota_1}{\iota_2}{\iota_2}{\iota_2}{\iota_2}{\iota_1}`r=(r_1, r_2)`{\columnII{x_1}{x_2}\times\columnII{y_1}{y_2}}]
\efig
$$

By the assumptions on $R$, the morphisms $(x_1, y_2)$, $(x_2, y_2)$ and $(x_2, y_1)$ from $S$ to $X \times Y$ all factor through $r$. It follows that there is a morphism $p$ as shown which makes the diagram commute.  By the property of the  factorisation, since $r$ is a monomorphism, there is a morphism $f: R' \rightarrow R$ which makes the following diagram commute:
$$
\bfig
\square(0,0)|alrb|<1000,500>[R'`R`2S \times 2S` X \times Y; f`r' = (r_1', r_2')`%
r=(r_1, r_2)`\columnII{x_1}{x_2}\times\columnII{y_1}{y_2}]
\efig
$$
By hypothesis, $\iota_1$ and $\iota_1$ are related by $R'$, so that the map $(\iota_1, \iota_1): S \rightarrow 2S \times 2S$ factors through $r'$. By commutativity of the above diagram, the map $(x_1, y_1): S \rightarrow X\times Y$ must then factor through $r$, as required.
\end{proof}

Note that Proposition \ref{bojan1} holds more generally for any category with (strong epi, mono)-factorizations and binary products and coproducts, where one replaces the (regular epi, mono)-factorization of the vertical morphism in (\ref{eq:fact1}) with its (strong epi, mono)-factorization.

\begin{corollary}\label{largest}
Let $\mathbb{C}$ be a regular category admitting binary coproducts. Let $\mathbb{D}$ be a full subcategory of $\mathbb{C}$ which is a Mal'tsev category and which is closed under regular quotients and binary coproducts in $\mathbb{C}$. Then $\mathbb{D}$ is contained in $\mathsf{Mal}(\mathbb{C})$.
\end{corollary}
\begin{proof}
Suppose that $\mathbb{D}$ is a full subcategory of $\mathbb{C}$ which is Mal'tsev and which is closed under binary coproducts and regular quotients in $\mathbb{C}$. Then for every object $S$ in $\mathbb{D}$, the objects $2S$ and $R'$ from diagram (\ref{eq:fact1}), and hence also the morphisms $r'_1$ and $r'_2$, are contained in $\mathbb{D}$. Since the morphisms $r'_1$ and $r'_2$ are jointly monic in $\mathbb{C}$, they are also jointly monic in $\mathbb{D}$ and thus represent a internal relation in $\mathbb{D}$. Since $\mathbb{D}$ is assumed to be Mal'tsev, the relation $\Hom(S,R')$ between sets must be difunctional. But then, since $\iota_1 R' \iota_2$, $\iota_2 R' \iota_2$ and $\iota_2 R' \iota_1$, we have that  $\iota_1 R' \iota_1$, so $S$ is a Mal'tsev object by Proposition~\ref{bojan1}.
\end{proof}

It is not clear in general if the full subcategory $\mathsf{Mal}(\mathbb{C})$ is itself a Mal'tsev category, since jointly monic pairs in $\mathsf{Mal}(\mathbb{C})$ may not be jointly monic as morphisms in $\mathbb{C}$. The following corollary gives a condition under which $\mathsf{Mal}(\mathbb{C})$ is indeed a Mal'tsev category.

\begin{corollary}\label{regepisagree}
Let $\mathbb{C}$ be regular category with binary coproducts. Consider the following conditions on $\mathbb{C}$:
\begin{itemize}
\item[(1)] every morphism in $\mathsf{Mal}(\mathbb{C})$ which is a regular epimorphism in $\mathbb{C}$ is also a regular epimorphism as a morphism in $\mathsf{Mal}(\mathbb{C})$;
\item[(2)] if a pair of morphisms $r_1: R \rightarrow X$ and $r_2: R \rightarrow Y$ are jointly monic (that is, an internal relation) in $\mathsf{Mal}(\mathbb{C})$ then they are also jointly monic in $\mathbb{C}$;
\item[(3)] $\mathsf{Mal}(\mathbb{C})$ is the largest full subcategory of $\mathbb{C}$ which is a Mal'tsev category and which is closed under regular quotients and binary coproducts in $\mathbb{C}$.
\end{itemize}
Then (1) $\Rightarrow$ (2) and (2) $\Rightarrow$ (3).
\end{corollary}
\begin{proof}
(1) $\Rightarrow$ (2):  Suppose $r_1: R \rightarrow X$ and $r_2: R \rightarrow Y$ are jointly monic in $\mathsf{Mal}(\mathbb{C})$. Consider the map $(r_1, r_2): R \rightarrow X\times Y$ in $\mathbb{C}$ and its (regular epi, mono)-factorization $(r_1, r_2) = me$ in $\mathbb{C}$. Since $e$ is a regular epi in $\mathbb{C}$, it exists as a morphism in $\mathsf{Mal}(\mathbb{C})$ where it is also a regular epi. Since $r_1$ and $r_2$ are jointly monic in  $\mathsf{Mal}(\mathbb{C})$, it is easy to check that $e$ must be a monomorphism in  $\mathsf{Mal}(\mathbb{C})$, so that in fact $e$ is an isomorphism. It follows that $(r_1, r_2)$ is a monomorphism in $\mathbb{C}$ as required.

(2) $\Rightarrow$ (3): Since internal relations in $\mathsf{Mal}(\mathbb{C})$ are internal relations in $\mathbb{C}$, $\mathsf{Mal}(\mathbb{C})$ is a Mal'tsev category, and the result follows from Corollary \ref{largest}.
\end{proof}

Conditions (1) and (2) in Corollary \ref{regepisagree} will turn out to hold for the categories we are interested in in the next two sections  ($\mathbf{Top}^\mathrm{op}$ and $\mathbf{Met}_\infty^\mathrm{op}$). However, they are not very natural to require of a general category $\mathbb{C}$. 

\begin{question}
Are there natural conditions on a general category $\mathbb{C}$ such that  $\mathsf{Mal}(\mathbb{C})$ is a Mal'tsev category?
\end{question}

It follows from the work in~\cite{BoJan08} that, for a finitely cocomplete regular category $\mathbb{C}$ and an object $S$ in $\mathbb{C}$, conditions (a) and (b) in Proposition~\ref{bojan1} are further equivalent to the following:
\begin{itemize}
\item[(c)] $S$ admits an \emph{approximate Mal'tsev co-operation} $\mu$ \emph{with approximation} $\alpha$ a regular epimorphism, i.e.~there exists an object $A$ and morphisms $\mu: A \rightarrow 3S$ and $\alpha: A \rightarrow S$, with $\alpha$ a regular epimorphism, such that the following diagram commutes:
\begin{equation}\label{eq:approx1}
\bfig
\square(0,0)|alrb|<1000,500>[A`3S`S`2S \times 2S; %
\mu_S`\alpha_S`\matrIIxIII{\iota_1}{\iota_2}{\iota_2}{\iota_2}{\iota_2}{\iota_1}` (\iota_1,\iota_1)]
\efig
\end{equation}
\end{itemize}
Indeed, if such a diagram (\ref{eq:approx1}) exists with $\alpha_S$ a regular epimorphism, then by the universal property of the (regular epi, mono)-factorization system, the morphism $(\iota_1, \iota_1)$ factors through the monomorphism $(r_1', r_2')$ in (\ref{eq:fact1}), which implies (b). Conversely, if (b) holds, there is a map $g: S \rightarrow R'$ such that $(r_1', r_2') \circ g = (\iota_1, \iota_1)$, and one can take $\alpha_S$ to be the pullback of the map $e$ in (\ref{eq:fact1}) along the map $g$. Since $\mathbb{C}$ was assumed to be regular, $\alpha_S$ is a regular epimorphism.

\section{Co-Mal'tsev objects in $\mathbf{Top}$}\label{r1}
It is easy to check that the regular monomorphisms in $\mathbf{Top}$ are precisely the embeddings of spaces. In particular, an embedding $f: A \rightarrow X$ is the coequalizer of the continuous maps $a$ and $b$ to $J = \{0,1\}$, the two element indiscrete space, where $a$ sends $f(A)$ to $0$ and its complement to $1$ and $b$ sends all of $X$ to $0$. It is easy to show, moreover, that topological embeddings are closed under pushouts in $\mathbf{Top}$; it follows that the category $\mathbf{Top}^\mathsf{op}$ is regular and finitely complete.

\begin{theorem}\label{main}
Let $S$ be an object in $\mathbf{Top}$, and let the following diagram in $\mathbf{Top}$ represent the (epi, regular mono)-factorization of  the vertical morphism:
$$
\bfig
\Dtriangle(0,0)|lba|<400,-300>[S^2+S^2`R'`S^3;%
\matrIIIxII{\pi_1}{\pi_2}{\pi_2}{\pi_2}{\pi_2}{\pi_1}`r'`\rowIII{k_1}{k_2}{k_3}]
\efig
$$
Then the following are equivalent:
\begin{itemize}
\item[(a)] $S$ is a co-Mal'tsev object;
\item[(b)] there is a (unique) morphism $f: R' \rightarrow S$ such that 
$$f \circ r' = \columnII{\pi_1}{\pi_1};$$
\item[(c)] for every open set $A$ in $S$, there is an open set $A'$ in $S^3$ such that 
$$
x \in A \Leftrightarrow (x,y,y) \in A' \Leftrightarrow (y,y,x) \in A'
$$
for all $(x,y) \in S^2$;
\item[(d)] $S$ is an $R_1$-space.
\end{itemize}
\end{theorem}
\begin{proof}
(a) $\Leftrightarrow$ (b) follows from the dual of Proposition~\ref{bojan1}.

(b) $\Leftrightarrow$ (c): Since $(k_1, k_2, k_3)$ is a regular monomorphism, i.e.~an embedding of spaces, $R'$ has underlying set
$$ \{(x,y,y) \mid (x, y) \in S^2 \} \cup \{(y,y,x) \mid (x, y) \in S^2\} $$
with the subspace topology induced by $S^3$. Let $f$ be the function from the underlying set of $R'$ to the underlying set of $S$ defined by $f(x,y,y) = x$ and $f(y,y,x) = x$. Condition (b) is then equivalent to $f$ being a continuous map from $R'$ to $S$, which is clearly equivalent to (c).

(c) $\Rightarrow$ (d): Let $x,y$ be two points in $S$ and let $A$ be an open set such that $x \in A$ and $y \notin A$. Then take $A'$ as in (c). Now $(x,y,y) \in A'$, so there exist open sets $U$, $V$, and $W$ in $S$ such that $(x,y,y) \in U\times V \times W \subseteq A'$. Moreover, $(x,y,y) \in U\times (V\cap W) \times (V \cap W) \subseteq A'$. Now suppose $z \in U\cap V \cap W$. Then $(z,z,y)$ is in $A'$ and thus $y$ must be in $A$, a contradiction. So $U\cap (V \cap W) = \varnothing$ and thus $U$ and $V \cap W$ are disjoint open sets such that $x \in U$ and $y \in V\cap W$.

(d) $\Rightarrow$ (c): Let $A$ be an open set in $S$. Let $(x,y)$ be a pair of points with $x \in A, y \notin A$. Then since $S$ is an $R_1$-space, there exist disjoint open sets $B_{(x,y)}, C_{(x,y)}$ such that $x \in B_{(x,y)} \subseteq A$ and $y \in C_{(x,y)}$. Now consider the family of all such pairs $(B_{(x,y)}, C_{(x,y)})$ indexed by pairs of points $(x,y)$ with $x \in A$ and $y \notin A$. Now it is easy to see that the desired set $A'$ may be chosen to be :

$$
A' = A^3 \cup \left(\bigcup_{x\in A, y \notin A}  B_{(x,y)} \times C_{(x,y)} \times C_{(x,y)} \right) \cup \left(\bigcup_{x\in A, y \notin A} C_{(x,y)} \times C_{(x,y)} \times B_{(x,y)}  \right)
$$
\end{proof}

We thus have the following corollaries of Corollary~\ref{reflect}, Corollary~\ref{regepisagree} and the remark on approximate Mal'tsev co-operations at the end of the previous section.

\begin{corollary}
Let $\mathbf{R}_1$ be the full subcategory of $\mathbf{Top}$ whose objects are the $R_1$-spaces. Then the dual of $\mathbf{R}_1$ is a finitely complete Mal'tsev category. Moreover, $\mathbf{R}_1$ is reflective in $\mathbf{Top}$ and is the largest full subcategory of $\mathbf{Top}$ whose dual is Mal'tsev and which is closed under binary products and regular subobjects (i.e.~subspaces) in $\mathbf{Top}$.
\end{corollary}
\begin{proof}
The only part which needs proving is that if a morphism $f$ in $\mathbf{R}_1$ is a regular monomorphism in $\mathbf{Top}$ then it is a regular monomorphism in $\mathbf{R}_1$, after which one can apply Corollary \ref{regepisagree}. This is easy to check given that the two element indiscrete space is in $\mathbf{R}_1$.
\end{proof}

The notion of \emph{approximate Mal'tsev operation} is dual to that of an approximate Mal'tsev co-operation.

\begin{corollary}\label{ops}
Let $X$ be an object of $\mathbf{Top}$. Then $X$ is an $R_1$-space if and only if $X$ admits an approximate Mal'tsev operation with approximation $\alpha$ a regular monomorphism (i.e. an embedding of spaces).
\end{corollary}

\section{Co-Mal'tsev objects in $\mathbf{Met}$}\label{ultra}
Let $\mathbf{Met}$ be the category whose objects are metric spaces and whose morphisms are all \emph{short maps} between metric spaces, i.e.~maps $f: X \rightarrow Y$ such that for any $x_1$, $x_2 \in X$, 
$$
d(f(x_1), f(x_2)) \leqslant d(x_1, x_2).
$$
This is, for example, the category of metric spaces implicit in Isbell's definition of injective metric space in \cite{Isbell64}. Note that short maps are always continuous with respect to the topology induced by the metric on each space, and that the isomorphisms in $\mathbf{Met}$ are precisely the global isometries. In this section we will prove that the co-Mal'tsev objects in this category are the ultrametric spaces. 

The category $\mathbf{Met}$ does not admit coproducts, so we will also want to consider the category $\mathbf{Met}_\infty$ whose objects are \emph{extended metric spaces} and whose maps are short maps between extended metric spaces. Recall that an extended metric space is a set equipped with a distance function which takes values in the extended reals $\mathbb{R} \cup \{\infty\}$ and which satisfies the axioms for a metric. In particular, every metric space can be viewed as an extended metric space. We now collect some elementary facts about $\mathbf{Met}_\infty$, leading eventually to Proposition \ref{metopisregular} below. The results are straightforward to prove, but we include proofs for the sake of completeness.

To construct colimits in $\mathbf{Met}_\infty$ we need to be able to take quotients of metric spaces by equivalence relations. This topic is classical (see for example~\cite{Burago, Gromov07}). Let $A$ be a metric space and $E$ an equivalence relation on $A$. Let $A_E$ be the set of equivalence classes under $E$ and define a distance function $d'$ on $A_E$ as follows:
$$
d'([a]_E, [b]_E) = \mathsf{inf}\{ \sum_{i=1}^{n} d(a_i, b_i) \mid a_1Ea,\ b_nEb,\ b_i E a_{i+1},\ n \in \mathbb{Z}_{+}  \}
$$
where $\mathbb{Z}_{+}$ is the set of positive integers. A sequence of pairs $(a_i, b_i)_{1 \leq i \leq n}$ satisfying $a_1Ea,\ b_nEb,\ b_i E a_{i+1}$ will usually be referred to as a \emph{chain} from $a$ to $b$. In general this defines a pseudometric on $A_E$, which may not be a metric (some distinct points may be distance $0$ apart). If $d'$ is a metric, then we define $\overline{A}_E$ to be $A_E$ with the metric $d'$; in such a case, we will call the equivalence relation $E$ \emph{well-behaved}. If $d'$ is not a metric, consider the equivalence relation $x \sim y \Leftrightarrow d'(x,y) = 0$ on $A_E$, and define $\overline{A}_E$ to be the set of equivalence classes under $\sim$ with the metric
$$
d_E([x]_{\sim}, [y]_{\sim} ) = d'(x, y)
$$
(note that this is well-defined). It is an easy exercise to check that $\overline{A}_E$ with the obvious quotient map is universal amongst all short maps with domain $A$ which are constant on equivalence classes under $E$.

Using this construction it is easy to define coequalizers in $\mathbf{Met}$ and $\mathbf{Met}_\infty$: for two maps $f,g: X \rightarrow Y$ simply take the quotient of $Y$ by the equivalence relation generated by the pairs $(f(x), g(x))_{x \in X}$. Coproducts in $\mathbf{Met}_\infty$ are easy to construct (but don't exist in $\mathbf{Met}$):  to form $X+Y$ simply take the disjoint union of the two spaces and declare the distance between any point in $X$ and any point in $Y$ to be infinite. It follows that $\mathbf{Met}_\infty$ is finitely cocomplete. It is easy to check that $\mathbf{Met}$ and $\mathbf{Met}_\infty$ also admit equalizers. Given a pair of objects $X$ and $Y$ in either $\mathbf{Met}$ or $\mathbf{Met}_\infty$, their product is given by the set $X \times Y$ equipped with the metric
$$
d((x_1, y_1), (x_2, y_2)) = \mathsf{max}(d(x_1, x_2), d(y_1, y_2)).
$$

\begin{lemma}\label{regmonosinmet}
A morphism $f: X \rightarrow Y$ in $\mathbf{Met}$ or $\mathbf{Met}_\infty$ is a regular monomorphism if and only if it is an isometric embedding with closed image.
\end{lemma}
\begin{proof}
The same proof will work for both categories. Suppose $f$ is the equalizer of a pair $a,b: Y \rightarrow Z$. Since $a$ and $b$ agree on $f(X)$, they also agree on the closure $\overline{f(X)}$. It follows that $f(X) = \overline{f(X)}$ so that the image of $f$ is closed, and it is easy to check that $f$ must be an isometric embedding.

Conversely, suppose $f: X \rightarrow Y$ is an isometric embedding with $f(X)$ closed. The case when $X$ is empty is easy to check, so assume $x_0 \in X$. Consider the quotient $Z$ of $Y$ by the equivalence relation 
$$
y \sim y' \Leftrightarrow \{y, y'\}\subseteq f(X) \text{ or } y = y', 
$$
which one checks is well-behaved because $f(X)$ is closed. It is now easy to check that $f$ is the equalizer of the quotient map $a: Y \rightarrow Z$ and the map $b$ which sends all of $Y$ to $[f(x_0)]$.
\end{proof}

\begin{proposition}\label{metopisregular}
The dual of the category $\mathbf{Met}_\infty$ is a finitely complete and finitely cocomplete regular category. 
\end{proposition}
\begin{proof}
Given a morphism $f$ in $\mathbf{Met}_\infty$, we have an (epi, regular mono)-factorization $f = me$ given by
$$ \bfig
\Vtriangle(0,0)|abb|/{->}`{->}`{<-}/<400, 400>[X`Y`\overline{f(X)};%
f`e`m]
\efig
$$
where the metric on $\overline{f(X)}$ is inherited from $Y$. From this it is easy to show that $\mathbf{Met}_\infty$ admits a (regular epi, mono)-factorization system, and we have already noted that it is finitely complete and cocomplete. 

It remains to show that regular monos are closed under pushouts. Let $m: X \rightarrow Y$ be a regular monomorphism, which we may suppose is an inclusion of a closed subspace $X \subseteq Y$, and let $f: X \rightarrow Z$ an arbitrary morphism. The pushout of of $m$ along $f$ is given by the quotient of $Y + Z$ by the equivalence relation $\sim$ generated by the pairs $(x, f(x))_{x \in X}$. Denote this space by $Q$ and let $m': Z \rightarrow Q$, $f': Y \rightarrow Q$ be the maps induced by inclusions. We claim that $m'$ is an isometric embedding. Let $z, z' \in Z$ and let $(a_i, b_i)_{1 \leq i \leq n}$ be a chain in $Y+Z$ from $z$ to $z'$ with respect to $\sim$. We want to show that
$$
\sum_{i=1}^{n} d_{Y + Z}(a_i, b_i)  \geq d_Z(z, z').
$$
We may assume that none of the distances $d_{Y + Z}(a_i, b_i)$ is infinite. If some $a_i$ or $b_i$ is in $Y \setminus X \subseteq Y + Z$, then there is a subchain $(a_i, b_i)_{k \leq i \leq k'}$ with $a_k \in X$, $b_{k'} \in X$ and all other elements in $Y \setminus X$. Since $\sim$ is trivial on $Y \setminus X$, we have $b_i = a_{i+1}$ for $k \leq i \leq k'-1$, so by the triangle inequality
$$
\sum_{i=k}^{k'} d_{Y + Z}(a_i, b_i) \geq d_Y(a_k, b_{k'}) = d_X(a_k, b_{k'}).
$$
Thus we may eliminate subchains in $Y\setminus X$ to reduce to the case when all the $a_i$ and $b_i$ are in $X \sqcup Z \subseteq Y+Z$.

If all the $a_i$ and $b_i$ are in $X \sqcup Z \subseteq Y+Z$, then applying the short map $f'': X + Z \rightarrow Z$ induced by $f$ and $1_Z$, we obtain a chain $(f''(a_i), f''(b_i))_{1 \leq i \leq n}$ with $f''(a_1) = z$, $f''(b_n) = z'$ and $f''(b_i) = f''(a_{i+1})$, and by the triangle inequality
$$
\sum_{i=1}^{n} d_{Y + Z}(a_i, b_i) \geq \sum_{i=1}^{n} d_{Z}(f''(a_i), f''(b_i))
\geq d_{Z}(z, z').
$$
This concludes the proof that $d_Z(z,z') = d_{Q}(m'(z),m'(z'))$ for all $z,z' \in Z$. It remains to show that $m'(Z)$ is closed. To see this, let $y \in Y \setminus X \subseteq Y+Z$. Since $\sim$ is trivial on $Y\setminus X$, any chain $(a_i, b_i)_{1 \leq i \leq n}$ from $y$ to $z \in Z \subseteq Y+Z$ must have some minimal $j$ with $a_j \in Y \setminus X$ and $b_j \in X$. It follows that 
$$
\sum_{i=1}^{n} d_{Y + Z}(a_i, b_i) \geq \sum_{i=1}^{j} d_{Y + Z}(a_i, b_i) \geq d_Y(y, b_j)
$$
so that $d_Q(f'(y), m'(z))$ is bounded below by $d_Y(y, X)$, which is positive since $X$ is closed. Thus $f'(y)$ has a neighbourhood which does not intersect $m'(Z)$ as required.
\end{proof}

\begin{theorem}\label{ultramain}
Let $S$ be an object of $\mathbf{Met}$ or $\mathbf{Met}_\infty$. Then $S$ is a co-Mal'tsev object in the respective category if and only if it is an ultrametric space.
\end{theorem}
\begin{proof}
Consider, in $\mathbf{Met}_\infty$, the dual picture to diagram (\ref{eq:fact1}), namely the diagram
$$
\bfig
\Dtriangle(0,0)|lba|<400,-400>[S^2+S^2`R'`S^3;\matrIIIxII{\pi_1}{\pi_2}{\pi_2}{\pi_2}{\pi_2}{\pi_1}`r'=\columnII{r'_1}{r'_2}`\rowIII{k_1}{k_2}{k_3}]
\efig
$$
where $(k_1, k_2, k_3) \circ r'$ represents the factorisation of the vertical morphism into an epi followed by a regular mono. It follows from the proof of Proposition~\ref{metopisregular} that the object $R'$ is the closure of the subspace
$$
T = \{(x, y, y) \mid (x, y) \in S^2\} \cup \{(x, x, y) \mid (x, y) \in S^2\}
\subseteq S^3.$$
One can easily check, however, that $T$ is itself closed, so that $R'$ is just the subspace $T$. Thus the internal co-relation $(R', r'_1, r'_2)$ is the co-relation
\begin{equation}\label{eq:corel}
\bfig
\morphism(-500,0)[S^2`R';r_1']
\morphism(500,0)<-500,0>[S^2`R';r_2']
\efig
\end{equation}
with $R' = T$ and where the maps $r'_1$ and $r'_2$ send $(x,y)$ to $(x, y, y)$ and $(x, x, y)$ respectively. We are now ready to prove the theorem.

$(\Rightarrow)$ If $S$ is a metric space, then the space $R'$ in (\ref{eq:corel}) is a metric space, as is $S^2$, so we can form the co-relation $R'$ in diagram (\ref{eq:corel}) in both $\mathbf{Met}$ and $\mathbf{Met}_\infty$. We have $\pi_1 R' \pi_2$, $\pi_2 R'  \pi_2$ and $\pi_2 R' \pi_1$, so if $S$ is a co-Mal'tsev object then there must exist a morphism $f: R' \rightarrow S$ such that $f \circ r'_1 = f \circ r'_2 = \pi_1$. The map $f$ is uniquely defined: it sends $(x,y,y)$ to $x$ and $(x, x, y)$ to $y$. Let $x,y,z$ be points in $S$. Then since $f$ is a short map, we have
\begin{align*}
d(x,z) & = d(f(x, y, y), f(y, y, z))  \\
& \leqslant d((x,y, y), (y,y,z)) \\
& = \mathsf{max}\{d(x, y), d(y, y) , d(y, z) \} \\
& = \mathsf{max}\{d(x, y), d(y, z) \} .\\
\end{align*}

$(\Leftarrow)$ Let $S$ an object of $\mathbf{Met}_\infty$ which is an ultrametric space, and let the internal co-relation $R'$ be as above. By the above results, the dual of $\mathbf{Met}_\infty$ satisfies the assumptions of Proposition~\ref{bojan1}, so it is enough to show the existence of a map $f: R' \rightarrow S$ such that $f \circ r'_1 = f \circ r'_2 = \pi_1$. Define $f$ to send $(x, y, y)$ to $x$ and $(x, x, y)$ to $y$. It remains to show that $f$ is a short map. We have
\begin{align*}
d(f(u,v,v), f(x, x, y)) & = d(u, y) \\
& \leqslant \mathsf{max}\{d(u, v), d(v, y)\}  \\
& \leqslant \mathsf{max}\{d(u, x), d(x, v), d(v, y) \} \\
& = d((u, v, v), (x, x, y)), \\
\end{align*}
and it follows easily that $f$ is a short map. Thus $S$ is a co-Mal'tsev object of $\mathbf{Met}_\infty$. 

Suppose now that $S$ is actually an object of $\mathbf{Met}$. Since $S$ is a co-Mal'tsev object in $\mathbf{Met}_\infty$, it is enough to check that a pair of jointly epic morphisms (that is, an internal co-relation) $a: A \rightarrow X,\ b: B \rightarrow X$ in $\mathbf{Met}$ remains jointly epic in $\mathbf{Met}_\infty$. Suppose $f,g: X \rightarrow Y$ are short maps such that $Y$ is an $\infty$-metric space, $fa = fb$ and $ga = gb$. If $A = B =\varnothing$, then clearly $X = \varnothing$, so we may exclude this case. Since $X$ is a metric space, the images of $f$ and $g$ each lie in a subspace of $Y$ which is a metric space. In fact, since $f$ and $g$ agree on at least one point, we may choose the subspaces to be the same. It follows that $f$ and $g$ each factor through a monomorphism $Y' \rightarrow Y$ with $Y'$ a metric space, and so $f= g$ as required.
\end{proof}

Let $\mathbf{UMet}$ and $\mathbf{UMet}_\infty$ be the full subcategories of $\mathbf{Met}$ and $\mathbf{Met}_\infty$ respectively whose objects are the ultrametric spaces.

\begin{lemma}\label{umetepis}
Let $f: X \rightarrow Y$ be a regular monomorphism in $\mathbf{Met}$ or $\mathbf{Met}_\infty$ where $X$ and $Y$ are ultrametric spaces. Then $f$ is a regular monomorphism in $\mathbf{UMet}$ or $\mathbf{UMet}_\infty$ respectively.
\end{lemma}
\begin{proof}
The same proof works for both categories. By Lemma~\ref{regmonosinmet}, $f$ is an isometric embedding with closed image. We may assume that $X$ is non-empty, with $x_0 \in X$. Consider the space $Z$ whose underlying set is the quotient of $Y$ by the equivalence relation 
$$
y \sim y' \Leftrightarrow \{y, y'\}\subseteq f(X) \text{ or } y = y', 
$$
and whose metric is given by 
$$
d_Z([a]_{\sim}, [b]_{\sim}) = \mathsf{inf}\{ \max \{ d(a_i, b_i) \mid 1 \leq i \leq n\} \mid a_1Ea,\ b_nEb,\ b_i E a_{i+1},\ n \in \mathbb{Z}_{+}  \}.
$$
It is easy to see that $d_Z$ satisfies $d_Z(z,z'') \leq \max(d_Z(z,z'), d_Z(z',z''))$, so we want to check that $d_Z(z,z') = 0 \implies z = z'$. Let $y \neq y'$ be in $Y$ such that $[y] \neq [y']$, and let $(a_i, b_i)_{1 \leq i \leq n}$ be a chain from $y$ to $y'$ in $Y$. If all the $a_i$ and $b_i$ are in $Y \setminus f(X)$, then
$$
\max \{ d_Y(a_i, b_i) \mid 1 \leq i \leq n\} \geq d_Y(y, y') > 0
$$ 
because $\sim$ is trivial on $Y\setminus X$ and $d_Y$ is an ultrametric. Thus we may restrict to chains where some $b_i$ is in $f(X)$. Let $j$ be the minimal index for which $b_j$ is in $f(X)$. We have 
$$\max \{ d_Y(a_i, b_i) \mid 1 \leq i \leq j\} \geq d_Y(y, b_j) \geq d(y, f(X)) > 0$$
since $f(X)$ is closed. Thus $d_Z([y], [y']) > 0$ as required. Finally, it is easy to check that $f$ is the coequalizer of the maps $a$ and $b$ where $a: Y \rightarrow Z$ is the quotient map and $b$ sends all of $Y$ to $[f(x_0)]$.
\end{proof}

\begin{corollary}
$\mathbf{UMet}_\infty$ (resp. $\mathbf{UMet}$) is the largest  full subcategory of $\mathbf{Met}_\infty$ (resp. $\mathbf{Met}$) whose dual is a Mal'tsev category and which is closed under products and regular subsobjects (i.e.~isometric embeddings of closed subspaces) in $\mathbf{Met}_\infty$ (resp. $\mathbf{Met}$) . 
\end{corollary}
\begin{proof}
The result about $\mathbf{UMet}_\infty$ follows from Corollary \ref{regepisagree} and Lemma \ref{umetepis}. Since $\mathbf{Met}^\mathrm{op}$ is not finitely complete, we need to make some elementary arguments to prove the result for $\mathbf{Met}$. If $\mathbb{D}$ is a full subcategory of $\mathbf{Met}^\mathrm{op}$ which is Mal'tsev and closed under coproducts and regular quotients in $\mathbf{Met}^\mathrm{op}$, then it is also closed under coproducts and regular quotients in $\mathbf{Met}_\infty^\mathrm{op}$, so by the above assertion, it is contained in $\mathbf{UMet}_\infty^\mathrm{op}$, and hence in $\mathbf{UMet}_\infty^\mathrm{op} \cap \mathbf{Met}^\mathrm{op} = \mathbf{UMet}^\mathrm{op}$. It remains to show that $\mathbf{UMet}^\mathrm{op}$ is itself a Mal'tsev category. If $a,b$ is a jointly monic pair in $\mathbf{UMet}^\mathrm{op}$, then it is also jointly monic in $\mathbf{UMet}_\infty^\mathrm{op}$ by similar arguments to the end of the proof of Theorem \ref{ultramain}. It is thus also jointly monic in $\mathbf{Met}_\infty^\mathrm{op}$ by Lemmas \ref{regepisagree} and \ref{umetepis}. Thus since every object in $\mathbf{UMet}^\mathrm{op}$ is a Mal'tsev object in $\mathbf{Met}_\infty^\mathrm{op}$ and every internal relation in $\mathbf{UMet}^\mathrm{op}$ is an internal relation in $\mathbf{Met}_\infty^\mathrm{op}$, $\mathbf{UMet}^\mathrm{op}$ is a Mal'tsev category.
\end{proof}

\section{Other term conditions}
It would be interesting to see if there are other connections of the form of Theorem~\ref{main} and Theorem~\ref{ultramain} between well-known conditions from universal algebra and well-known conditions from topology and geometry. As remarked by the authors of~\cite{BoJan08}, it is straightforward to adapt Proposition~\ref{bojan1} by replacing the underlying notion of Mal'tsev category with another \emph{category with closed relations} in the sense of~\cite{JanClosed}. One of the examples mentioned in~\cite{BoJan08} is that of a \emph{subtractive category}~\cite{JanSub}. A variety of universal algebras is \emph{subtractive} in the sense of Ursini~\cite{Ursini} if its theory contains a constant $0$ and a binary term $s$ satisfying the term equations $s(x,x) = 0$ and $s(x, 0) = x$. A subtractive category can be defined as a pointed category $\mathbb{C}$ with finite limits such that every internal relation $R$ satisfies the following condition (see~\cite{JanClosed}):
\begin{equation} \label{eq:punc}
 xRx \wedge xR0  \Rightarrow 0Rx 
 \end{equation}
Note that $0$ denotes the zero morphism in the above condition. Consider now the following condition on an object $S$ in a pointed category $\mathbb{C}$:
\begin{itemize}
\item[$\mathsf{(S)}$] for any relation $(R, r_1, r_2)$ from an object $X$ to $X$, the following relation satisfies condition (\ref{eq:punc}) above (where $0$ is the zero morphism from $S$ to $X$):
$$ \bfig
\morphism(0,0)<-1000, 0>[\Hom(S, R)`\Hom(S, X) ;\Hom(S, r_1)] 
\morphism(0,0)<1000, 0>[\Hom(S, R)`\Hom(S, X);\Hom(S, r_2)]
\efig $$
\end{itemize}

Theorem~\ref{sub1} was originally proved by Z. Janelidze~\cite{ZJanUnpub}, in a form involving the analogue of approximate Mal'tsev operations for the subtractive case, and it served as the original inspiration for this paper. A sketch of the proof is given here.

\begin{theorem}\label{sub1}
Let $S$ be an object of the category $\mathbf{Top}_\ast$ of pointed topological spaces. Then $S$ satisfies $(\mathsf{S})$ as an object of the dual category $\mathbf{Top}_\ast^\mathsf{op}$ if and only if it satisfies the following condition, where $0$ is the base point of $S$:
\begin{itemize}
\item[$(\mathsf{S}')$] if $A$ is any open set and $x$ a point such that either $x\in A\ \wedge \ 0 \notin A$ or $x \notin A\ \wedge \ 0 \in A$, then there exists disjoint open sets $B$ and $C$ in $S$ such that $0 \in B$ and $x \in C$.
\end{itemize}
\end{theorem}
\begin{proof}
It follows from arguments similar to the proof of Proposition~\ref{bojan1} that for a pointed category $\mathbb{C}$ with binary products and coproducts and a (strong epi, mono)-factorization system, an object $S$ in $\mathbb{C}$ satisfies $(\mathsf{S})$ if and only if $0 R' 1_S$ , where $(R', r'_1, r'_2)$ is the internal relation from $S$ to $S$ appearing in the (strong epi, mono)-factorisation of the vertical morphism in the following diagram:
\begin{equation}\label{eq:fact2}
\bfig
\Dtriangle(0,0)|aab|<400,300>[2S`R'`S\times S;%
\matrIIxII{1_S}{1_S}{1_S}{0}`e`r'=(r'_1,r'_2)]
\efig 
\end{equation}
Consider the dual diagram (\ref{eq:fact2}) in the category of pointed topological spaces (which has products and which admits (epi, regular mono)-factorizations): 
$$
\bfig
\Dtriangle(0,0)|lba|<400,-300>[2S`R'`S^2;%
\matrIIxII{1_S}{1_S}{1_S}{0}`r'`\rowII{k_1}{k_2}]
\efig
$$
We see that the space $R'$ is the subspace of $S^2$ given by
$$
R' = \{(x,x) \mid x \in S \} \cup  \{(x, 0) \mid x \in S \}
$$
where $0$ is the base point of $S$. We conclude that the pointed topological space $S$  satisfies $(\mathsf{S})$ as an object of the dual category $\mathbf{Top}_\ast^\mathsf{op}$ if and only if the set map $g: R' \rightarrow S$ defined by $g(x,x) = 0$ and $g(x,0) = x$ is continuous. This is to say that for every open set $A \subseteq S$, there is an open set $A'$ in $S^2$ such that $A \times \{0\} = A' \cap R'$ if $0 \notin A$ and $A \times \{0\} \cup \Delta = A' \cap R'$ if $0 \in A$. It remains to check that this is equivalent to $\mathsf{(S')}$, which is not hard to show.  
\end{proof}

\section*{Acknowledgements}
This paper is based on research which the author carried out while he was a Masters student under the supervision of Z.~Janelidze. The author would like to thank the anonymous referee for their valuable suggestions and for identifying an important mistake in the first draft of this paper.

\end{document}